\documentclass[12pt]{elsarticle}

\usepackage{amsthm,float,amsmath,
amstext,amssymb,amsfonts,bbold,extsizes}
\usepackage{graphicx}
\usepackage[arrow,matrix,curve]{xy}

\usepackage[active]{srcltx}

 \newtheorem{theorem}{Theorem}
 \newtheorem{lemma}[theorem]{Lemma}
 \newtheorem{proposition}[theorem]{Proposition}
\theoremstyle{remark}
\newtheorem{remark}[theorem]{Remark}

\DeclareMathOperator{\diag}{diag}

\begin{document}

\sloppy

\title{Kleiner's theorem for unitary representations of posets}

\author{Yurii Samoilenko}
\address{Institute of Mathematics, Tereschenkivska 3, Kyiv, Ukraine.}
\ead{yurii.sam@imath.kiev.ua}

\author{Kostyantyn Yusenko\corref{cor}}
\address{Department of Mathematics, University of S\~ao Paulo, Brazil.}
\ead{kay.math@gmail.com}
\cortext[cor]{Corresponding author}

\begin{abstract}
A subspace representation of a poset $\mathcal S=\{s_1,\ldots,s_t\}$ is given by a system $(V;V_1,\ldots,V_t)$ consisting of a vector space $V$ and its subspaces $V_i$ such that  $V_i\subseteq V_j$ if $s_i \prec s_j$. For each real-valued vector $\chi=(\chi_1,\ldots,\chi_t)$  with positive components, we define a unitary $\chi$-representation of $\mathcal S$ as a system $(U;U_1,\ldots,U_t)$ 
that consists of a unitary space $U$ and its subspaces $U_i$ such that $U_i\subseteq U_j$ if $s_i\prec s_j$
and satisfies $\chi_1 P_1+\dots+\chi_t P_t= \mathbb 1$, in which $P_i$ is the orthogonal projection  onto $U_i$.

We prove that $\mathcal S$  has a finite number of unitarily nonequivalent indecomposable $\chi$-representations for each weight $\chi$ if and only if $\mathcal S$ has a finite number of nonequivalent indecomposable subspace representations; that is, if and only if $\mathcal S$ contains any of Kleiner's critical posets.
\end{abstract}


\maketitle

\section{Introduction}

Kleiner \cite{Kleiner1} described all partially ordered sets (\textit{posets}) with finite number of nonequivalent indecomposable representations. We extend his description to unitary representations of posets.

The notion of poset representations was introduced by Nazarova and Roiter \cite{NazarovaRoiter} (see also \cite{GabRoi,Simson1992}). A \textit{matrix representation} of a finite poset $\mathcal S=\{s_1,\ldots,s_t\}$ over a field $\mathbb F$ is a block matrix $\mathcal A=[A_1 | \ldots | A_t]$ over $\mathbb F$. Two representations
$\mathcal A=[A_1 | \ldots | A_t]$ and $\mathcal B=[B_1 | \ldots | B_t]$ are \textit{equivalent}  if $\mathcal A$ can be reduced to $\mathcal B$ by elementary row transformations, elementary column
transformations within $A_i$, and additions of
linear combinations of columns of $A_i$ to columns of $A_j$ if $s_i \prec s_j$.
The \textit{direct sum} of $\mathcal A$ and $\mathcal B$ is the representation
\[	
	\cal A\oplus \cal B:=
	\left[ \begin{array}{c c|c c| c | c c}
	         A_1 & 0  & A_2 & 0 & \ldots & A_t & 0 \\
	         0 & B_1 & 0 & B_2& \ldots & 0 & B_t
	        \end{array}
	\right].
\]
A representation is called \textit{indecomposable} if it is not equivalent to a direct sum of two representations. It is sufficient to classify only indecomposable representations since each representation is equivalent to a
direct sum of indecomposable representations, uniquely determined up to isomorphism of
summands.

Kleiner \cite{Kleiner1} (see also \cite[Theorem 5.1]{GabRoi} and \cite[Theorem 10.1]{Simson1992}) proved that a poset $\mathcal S$
has only a finite number of nonequivalent indecomposable  representations if and only if it does not contain a full poset whose Hasse diagram is one of the form
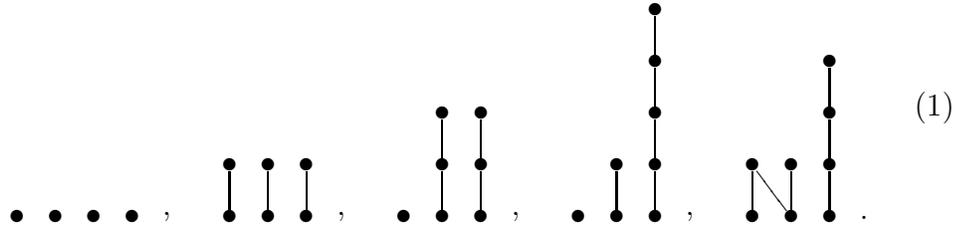
\begin{equation} \label{critial_posets}
\begin{split}
\entrymodifiers={[o]} \xymatrix @C=0.3cm @R=0.5cm
{
  & & & &  & &  &  &  &  &  & &  &  &  &    \bullet \ar@{-}[d] & & & &
\\
  & & & &  & &  &  &  &  & &  &  &  & & \bullet \ar@{-}[d] &&   &  & \bullet \ar@{-}[d]
\\
  & & & & &  &  &  &  & & \bullet \ar@{-}[d] & \bullet \ar@{-}[d] &  & & & \bullet \ar@{-}[d] & & &  & \bullet \ar@{-}[d]
  \\
  & & & & & \bullet \ar@{-}[d] & \bullet \ar@{-}[d] & \bullet \ar@{-}[d] & &  & \bullet \ar@{-}[d] & \bullet \ar@{-}[d] & & & \bullet \ar@{-}[d] &  \bullet \ar@{-}[d] & & \bullet \ar@{-}[d] \ar@{-}[rd] & \bullet \ar@{-}[d] & \bullet \ar@{-}[d]\\
 \bullet & \bullet & \bullet & \bullet &, \resizebox{0.3cm}{!}{ } & \bullet  & \bullet & \bullet & , \resizebox{0.3cm}{!}{ }  & \bullet &\bullet & \bullet & , \resizebox{0.3cm}{!}{ }  & \bullet & \bullet & \bullet & , \resizebox{0.3cm}{!}{ }  & \bullet & \bullet & \bullet&. }
\end{split}
\end{equation}

An equivalent definition of poset representations can be given in terms of subspaces. A \textit{subspace representation}  of $\mathcal S=\{s_1,\ldots,s_t\}$ is a tuple  $\mathcal V=(V;V_1,\ldots,V_t)$, in which $V$ is a vector space over $\mathbb F$ and $V_1,\ldots,V_t$ are its subspaces such that $V_i \subseteq V_j$ if $s_i \prec s_j$ (i.e., each representation is a homomorphism from $\mathcal S$ to the poset of all subspaces of  $V$). Two subspace representations $\mathcal V=(V;V_1,\ldots,V_t)$ and $\mathcal W=(W;W_1,\ldots,W_t)$ are \textit{equivalent}
if there exists a linear bijection
$g: V \to W$ such that $g(V_i)=W_i$ for all
$i$. For each subspace representation $\mathcal V=(V;V_1,\ldots,V_t)$, one can construct a matrix representation $\mathcal A=[A_1|\ldots|A_t]$ in such a way that (i) for each $i$ the columns of all $A_j$ with $s_j\preceq s_i$ generate the subspace $V_i$; and (ii) two subspace representations are equivalent if and only if the corresponding matrix representations are equivalent; see \cite[Chapter 3]{Simson1992}.

\medskip

From now on, all representations that we consider are over the field $\mathbb C$ of complex numbers. By a \textit{unitary representation of dimension} $d$, we mean a subspace representation $\mathcal U=(U;U_1,\ldots,U_t)$ in which $U$ is a unitary space of dimension $d$.  Two unitary representations $\mathcal U=(U;U_1,\ldots,U_t)$ and $\mathcal V=(V;V_1,\ldots,V_t)$  of a poset $\mathcal S$ are \textit{unitarily equivalent} if there exists a unitary bijection $\varphi: U \rightarrow V$ such that $\varphi(U_i)=V_i$ for all $i$.
The \textit{orthogonal sum} of unitary representations  $\mathcal U$ and $\mathcal V$ is the unitary representation
\[
\mathcal U\perp \mathcal V:=(U\perp V;U_1\perp V_1,\ldots,U_t\perp V_t),
\]
in which $U\perp V$ denotes the orthogonal sum of $U$ and $V$.
A unitary representation is called \emph{orthogonaly indecomposable} if it is not equivalent  to an orthogonal sum of two unitary representations.

Note that the problem of classifying unitary representations is hopeless even for the poset $\mathcal S =\{s_1,s_2,s_3 \ |\ s_1\prec s_2\}$ since by \cite[Theorem 4]{KruglyaSamoilenko2} it contains the problem of classifying an operator on a unitary space, and hence it contains the problem of  classifying any system of operators on unitary spaces \cite{KruglyaSamoilenko2,Sergeichuk}. The classification becomes possible for a broader class of posets if we impose additional conditions on unitary representations.

We denote the orthogonal projection onto
a subspace $M \subset U$ by $P_M$ and  the set of positive real numbers by $\mathbb R_+$.  We say that a unitary representation $\mathcal U=(U;U_1,\ldots,U_t)$ is a \textit{representation of weight} $\chi=(\chi_1,\ldots,\chi_t)\in \mathbb R_+^t$  (or \textit{$\chi$-representation}) if
\begin{equation} \label{orthoscalarEq}
	\chi_1 P_{U_1}+\dots+\chi_t P_{U_t}= \mathbb 1;
\end{equation}
such relations appear in many areas of mathematics, see for example \cite{AlbeverioOstrSamoilenko,Klyachko, KruglyaNazarovaRoiter,Totaro1994,Wu} and references therein.

Our goal is to prove that Kleiner's theorem holds for $\chi$-representations too:

\begin{theorem} \label {maintheoremFin}
The following conditions are equivalent for each finite poset $\mathcal S$ with $t$ elements:
\begin{enumerate}
\item[{\rm(i)}] For each $\chi \in \mathbb R_+^t$, $\mathcal S$ has only a finite number of indecomposable unitarily nonequivalent $\chi$-representations;
\item[{\rm(ii)}] For each $\chi \in \mathbb R_+^t$ and $d \in \mathbb N$, $\mathcal S$ has only a finite number of  indecomposable unitarily nonequivalent $\chi$-representations of dimension $d$;
\item[{\rm(iii)}] $\mathcal S$ does not contain a full poset whose Hasse diagram is one of the form \eqref{critial_posets}.
\end{enumerate}
\end{theorem}

\section{Preliminaries}

In what follows we suppose that the elements of a poset $\mathcal S$ are numbered from $1$ to $|\mathcal S|$. A poset is called \emph{primitive} and is denoted by $(t_1,\ldots,t_s)$ if it is the disjoint (cardinal) sum of linearly ordered sets of orders $t_i$.
The diagrams \eqref{critial_posets} and corresponding posets are called \textit{critical}. The poset which corresponds to the last diagram in the list   \eqref{critial_posets} is denoted by $(N,4)$.
To simplify the notation we denote a subspace representation $(V;V_1,\ldots,V_t)$ of $\mathcal S$ by $(V;V_i)_{i\in \mathcal S}$. The similar notation will be used for unitary representations and weights.

A subspace representation $\mathcal V=(V;V_i)_{i\in \mathcal S}$ is called \emph{schurian}
if all its endomorphisms are trivial; that is,
the ring ${\rm{End}}(\mathcal V):=\{ g \in M_{\dim  V}(\mathbb F)\ |\ g(V_i)\subseteq V_i,\ i \in \mathcal S\}$ is isomorphic to $\mathbb F$. Any schurian representation is indecomposable.

Any unitary representation $\mathcal U=(U;U_i)_{i \in \mathcal P}$ can be viewed
as a subspace representation;
the forgetful map is denoted by $F$.
If $\mathcal U$ is an indecomposable $\chi$-representation, then $F(\mathcal U)$ is schurian (see \cite[Theorem 1]{KruglyaNazarovaRoiter}).

\begin{lemma} \label{lemDiag}
	Let $P_i,Q_i \in M_n(\mathbb C)$, $i=1,\ldots,m$ be orthogonal projections such that
	\begin{equation} \label{lemAs}
	\chi_1 P_1+\dots+\chi_m P_m=\chi_1 Q_1+\dots+\chi_m Q_m
	\end{equation} for
	$(\chi_1,\ldots,\chi_m)$ with positive real $\chi_i$.
Let there exist a diagonal matrix $D=\diag(r_1,\ldots,r_n)$ with positive components such that $P_i D Q_i=D Q_i$ for all $i$. Then $r_1=\dots=r_n$ and $P_i=Q_i$ for all $i$.
\end{lemma}

\begin{proof}
	Write $P_i=[p^{(i)}_{k,l}]$, $Q_i=[q^{(i)}_{k,l}]$, $P_iDQ_i=[t^{(i)}_{k,l}]$, where $t^{(i)}_{k,l}=\sum_{j=1}^{n} r_j p^{(i)}_{k,j} \overline{q^{(i)}_{l,j}}$.
	Without losing generality, we may assume that $r_1=\max \{r_1,\ldots,r_n\}$. Since
	$D\sum_{i=1}^{m} \chi_i Q_i=\sum_{i=1}^{m} \chi_i P_i D Q_i$, we have
	\begin{equation} \label{firstIn}
	r_1\sum_{i=1}^{m} \chi_i q^{(i)}_{1,1}=\sum_{i=1}^{m} \chi_i \sum_{j=1}^{n} r_j p^{(i)}_{1,j} \overline{q^{(i)}_{1,j}}\leq r_1 \left |\sum_{i=1}^{m} \chi_i \sum_{j=1}^{n} p^{(i)}_{1,j} \overline{q^{(i)}_{1,j}}\right |.
	\end{equation}
For
	\begin{align*}
		 x&:=[\sqrt{\chi_1}p^{(1)}_{1,1},\ldots,\sqrt{\chi_1}p^{(1)}_{1,n},\ldots,\sqrt{\chi_m}p^{(m)}_{1,1},\ldots,\sqrt{\chi_m}p^{(m)}_{1,n}]^T \in \mathbb C^{nm},\\
		 y&:=[\sqrt{\chi_1}q^{(1)}_{1,1},\ldots,\sqrt{\chi_1}q^{(1)}_{1,n},\ldots,\sqrt{\chi_m}q^{(m)}_{1,1},\ldots,\sqrt{\chi_m}q^{(m)}_{1,n}]^T \in \mathbb C^{nm},
	\end{align*}
we have \[(x,y)=\sum_{i=1}^{nm}x_i\overline y_i=\sum_{i=1}^{m} \chi_i \sum_{j=1}^{n} p^{(i)}_{1,j} \overline{q^{(i)}_{1,j}}.\]
	Note that
	\[
	 |x|^2=\sum_{i=1}^{nm}x_i\overline x_i=\sum_{i=1}^m\chi_i\sum_{j=1}^n p^{(i)}_{1,j}\overline{p^{(i)}_{1,j}}=\sum_{i=1}^{m} \chi_i p^{(i)}_{1,1}.\] Similarly, $|y|^2=\sum_{i=1}^{m} \chi_i q^{(i)}_{1,1}$. By \eqref{lemAs}, we have  $\sum_{i=1}^{m} \chi_i p^{(i)}_{1,1}=\sum_{i=1}^{m} \chi_i q^{(i)}_{1,1}$, hence $|x|=|y|$. By the Cauchy-Schwartz inequality,
	\[
		\left |\sum_{i=1}^{m} \chi_i \sum_{j=1}^{n} p^{(i)}_{1,j} \overline{q^{(i)}_{1,j} } \right |\leq 	 \sqrt{\sum_{i=1}^{m} \chi_i p^{(i)}_{1,1}}\sqrt{\sum_{i=1}^{m} \chi_i q^{(i)}_{1,1}}=\sum_{i=1}^{m} \chi_i q^{(i)}_{1,1}.
	\]
	Comparing this inequality with (\ref{firstIn}), we have $r_1=\dots=r_n$ and $P_i=Q_i$ for all
	$i$.
 \end{proof}

\begin{theorem} \label{unitvslin}

Two $\chi$-representations $\mathcal U=(U;U_i)_{i\in \mathcal S}$ and $\mathcal U^\prime=(U^\prime;U_i^\prime)_{i\in \mathcal S}$  are unitarily equivalent if and only if
the corresponding subspace representations $F(\mathcal U)$ and $F(\mathcal U^\prime)$ are equivalent.
\end{theorem}

\begin{proof}
    If $\mathcal U$ is unitarily equivalent to $\mathcal U^\prime$, then $F(\mathcal U)$ is equivalent to $F(\mathcal  U^\prime)$. Let us prove the converse statement. $F(\mathcal U)$ is equivalent to $F(\mathcal U^\prime)$ if and only if there exists an invertible $g:U \to U^\prime$ such that
    \[
    	g^{-1} P_{U_i^\prime} g P_{U_i}=P_{U_i}, \
    	g P_{U_i} g^{-1} P_{U_i^\prime}=P_{U_i^\prime}, \quad i\in \mathcal S.
    \]
    Let $g=\varphi \psi D \psi^*$ be the polar decomposition of $g$, where $\varphi:U\to U^\prime$ and $\psi: U \to U$ are unitary maps and $D$ is a positively defined diagonal operator. Then
    \[
    	(\psi D^{-1} \psi^* \varphi^* )P_{U_i^\prime} (\varphi \psi D \psi^* )P_{U_i}=P_{U_i}, \
    	\quad i\in \mathcal S.
    \]
    Hence $(\psi^* \varphi^* P_{U_i^\prime} \varphi \psi) D (\psi^* P_{U_i} \psi)=D (\psi^* P_{U_i} \psi)$ for all  $i$. Since $\mathcal U$ and $\mathcal U^\prime$ are $\chi$-representations,
    \[
       \sum_{i\in \mathcal S} \chi_i (\psi^* \varphi^* P_{U_i^\prime} \varphi \psi) = \chi_0 I, \quad
    \sum_{i\in \mathcal S} \chi_i (\psi^* P_{U_i} \psi) = \chi_0 I.
    \]
    Lemma \ref{lemDiag} ensures that $\psi^* \varphi^* P_{U_i^\prime} \varphi \psi=\psi^* P_{U_i} \psi$ for all $i$. Therefore, $\varphi^* P_{U_i^\prime} \varphi =P_{U_i} $ for all $i$, and so $\mathcal U$ is unitarily equivalent to $\mathcal U^\prime$.
\end{proof}

\begin{remark} 
By similar argumentation, one can show that $\chi$-repre\-sentation $\mathcal U$ is orthogonally indecomposable if and only $F(\mathcal U)$ is indecomposable. The connection between usual and orthoscalar representations of quivers was established \cite[Theorem 1]{KruglyaNazarovaRoiter} in the same way as in Theorem \ref{unitvslin}.
\end{remark}

A representation $\mathcal V=(V;V_i)_{i\in \mathcal S}$ of weight $\chi=(\chi_i)_{i \in \mathcal S}$ is called $\chi$-\textit{stable} if
$
	\sum_{i \in \mathcal S} \chi_i \dim V_i=\dim V
$
and 
\[
	\sum_{i \in \mathcal S} \chi_i \dim(V_i \cap M)<\dim M
\]
for any proper subspace $0\neq M\subset V$.
\begin{lemma} \label{stab1}
	If \ $\mathcal U=(U;U_i)_{i \in \mathcal S}$ is an indecomposable $\chi$-representation, then $F(\mathcal U)$ is $\chi$-stable.
\end{lemma}
\begin{proof}
Equating the traces of both sides in \eqref{orthoscalarEq},   we obtain $\sum_{i \in \mathcal S} \chi_i \dim U_i=\dim U.$
 If $M$ is any proper subspace of $U$, then $
    \sum_{i \in \mathcal S}\chi_i P_{U_i}P_M=P_M$.
Equating the traces of both sides in the last equality, we get
\[
    \sum_{i \in \mathcal S}\chi_i {\rm{tr}}(P_{U_i}P_M)=\dim M.
\]
By \cite[Theorem 2]{hal},
${\rm{tr}}(P_{M_1\cap M_2}) \leq {\rm{tr}}(P_{M_1}P_{M_2})$  for each two subspaces $M_1$ and $M_2$, and so
\[
  \sum_{i \in \mathcal S}\chi_i {\rm{tr}}(P_{U_i\cap M})\leq\sum_{i \in \mathcal S}\chi_i {\rm{tr}}(P_{U_i}P_M)=\dim M.
\]
It remains to prove that the last inequality is strict. Indeed, assume
that  ${\rm{tr}}(P_{U_i\cap M})={\rm{tr}}(P_{U_i}P_M)$ for all $i$. Then
each $P_{U_i}$ commutes with $P_M$. Hence the subspace
$M$ is invariant with respect to the projections $P_{U_i}$ and the representation $\mathcal U$ is decomposable. This contradicts the assumption.
\end{proof}

The converse statement to Lemma \ref{stab1} also holds: if a representation $\mathcal V=(V;V_i)_{i\in \mathcal S}$ is $\chi$-stable, then one can choose a scalar product in $V$ in such a way that $\mathcal V$ becomes a $\chi$-representation; see \cite[Theorem 3.5]{Hu2004}. Using results from \cite{Hu2004, Klyachko, Totaro1994}, one can prove the following theorem.
\begin{theorem} \label{stability-un}
	An indecomposable unitary representation $\mathcal U$ is a $\chi$-repre\-sentation if and only if the corresponding subspace representation $F(\mathcal U)$ is $\chi$-stable.
\end{theorem}

\section{Proof of Theorem \ref{maintheoremFin}}

The implication ${\rm(i)}\Rightarrow {\rm(ii)}$ is trivial.

\medskip

${\rm(iii)}\Rightarrow {\rm(i)}$. Assume that the Hasse diagram of $\mathcal S$ does not contain
any of critical diagrams \eqref{critial_posets}. If $\mathcal S$ has an infinite number of indecomposable unitarily nonequivalent $\chi$-representations for some weight $\chi$, then by Theorem \ref{unitvslin} it has an infinite number of nonequivalent indecomposable subspace representations. By Kleiner's theorem, $\mathcal S$ contains a critical diagram; a contradiction.

\medskip

${\rm(ii)}\Rightarrow {\rm(iii)}$. We say that a poset $\mathcal S$ is \textit{unitary representation-infinite} if there exist $d\in \mathbb N$ and $\chi^{\mathcal S}\in \mathbb R^{|S|}_+$ such that $\mathcal S$ has an infinite number of indecomposable unitarily nonequivalent $\chi^{\mathcal S}$-representations of dimension $d$. Our aim is to prove that critical posets are unitary representation-infinite. 

One can  show that critical primitive posets are unitary representation-infinite using \cite{AlbeverioOstrSamoilenko,KruglyakLivinskyi, OstrovskyiSamoilenko}.
Namely, there exists a correspondence between the $\chi$-representations of a given poset $\mathcal S$ and the representations of a certain $*$-algebra $\mathcal A_{\Gamma,\omega}$ associated with a star-shaped graph $\Gamma$, which is determined by the Hasse diagram of $\mathcal S$, and the parameter $\omega$ is determined by the weight $\chi$. If $\Gamma$ is  an extended Dynkin graph (which corresponds to some primitive critical $\mathcal S$), then one can choose the parameter $\omega$ such that $\mathcal A_{\Gamma,\omega}$ has an infinite number of unitarily nonequivalent irreducible representations. The complete description of such representations was given in \cite{AlbeverioOstrSamoilenko, KruglyakLivinskyi,OstrovskyiSamoilenko} (see also Remark \ref{RemFinale}). But we use another method that handles both primitive and non-primitive cases.

Denote by $e^{(n)}_i$ (or $e_i$ if no confusion can arise) the $n$-dimensional vector in which the $i$-th coordinate is $1$ and the others are $0$. Denote by $e_{i_1\ldots i_k}$ the vector $e_{i_1}+\dots+e_{i_k}$ and by $\langle
x_1,\ldots,x_m \rangle$ the vector space spanned by
$x_1,\ldots,x_m \in \mathbb C^n$.

For each critical poset $\mathcal S$, we define a family of its subspace representations ${\cal V}_{\lambda }(\mathcal S)$ that depend on a complex parameter
$\lambda \in \mathbb C$.
\begin{itemize}
  \item If ${\cal S}=(1,1,1,1)$, then ${\cal V}_{\lambda }(\cal S)$ consists of the space $\mathbb C^2$ and its subspaces
\[\xymatrix @-1pc {\langle e_1 \rangle & \langle e_2 \rangle & \langle e_1+e_2 \rangle & {\langle e_1+\lambda e_2 \rangle}}\]

\item If ${\cal S}=(2,2,2)$, then ${\cal V}_{\lambda }(\cal S)$ consists of the space $\mathbb C^3$ and its subspaces
\[
\xymatrix @-1pc { {\langle e_{123},e_1+\lambda e_3 \rangle}&\langle e_1,e_2 \rangle&\langle e_2,e_3 \rangle\\ \langle e_{123} \rangle\ar@{->}[u]&\langle e_1 \rangle\ar@{->}[u] &\langle e_3 \rangle\ar@{->}[u]}\]

\item If ${\cal S}=(1,3,3)$, then ${\cal V}_{\lambda }(\cal S)$ consists of the space $\mathbb C^4$ and its subspaces
\[
\xymatrix @-1pc {& \langle e_1,e_4,e_2+\lambda e_3 \rangle& \langle e_1,e_2,e_3 \rangle\\
  &\langle e_1,e_4 \rangle\ar@{->}[u]&\langle e_2,e_3 \rangle\ar@{->}[u]\\
  \langle e_{123},e_{24} \rangle&\langle e_4 \rangle\ar@{->}[u]&\langle e_3 \rangle \ar@{->}[u]}
  \]

\item If ${\cal S}=(1,2,5)$, then ${\cal V}_{\lambda }(\cal S)$ consists of the space $\mathbb C^6$ and its subspaces
\[
\xymatrix @-1pc {& &\langle e_1,e_2,e_3,e_4,e_5+\lambda e_6 \rangle\\
                  &&\langle e_1,e_2,e_3,e_4 \rangle\ar@{->}[u]\\
                  &&\langle e_2,e_3,e_4 \rangle\ar@{->}[u]\\
                  &\langle e_1,e_2,e_5,e_6 \rangle&\langle e_3,e_4 \rangle\ar@{->}[u]\\
                  \langle e_{123},e_{245},e_{16} \rangle&\langle e_5,e_6 \rangle\ar@{->}[u]&\langle e_4 \rangle\ar@{->}[u]\\}
                  \]

\item If ${\cal S}=(N,4)$, then ${\cal V}_{\lambda }(\cal S)$ consists of the space $\mathbb C^5$ and its subspaces
\[\xymatrix @-1pc {& &\langle e_1,e_2,e_3,e_4 \rangle\\
                 &&\langle e_2,e_3,e_4 \rangle\ar@{->}[u]\\
                  \langle e_{235},e_{134},e_5,e_3+\lambda e_4 \rangle&\langle e_1,e_2,e_5 \rangle&\langle e_3,e_4 \rangle\ar@{->}[u]\\
                  \langle e_{235},e_{134} \rangle\ar@{->}[u]&\langle e_5 \rangle\ar@{->}[ul]\ar@{->}[u]&\langle e_4 \rangle\ar@{->}[u]}\]
\end{itemize}
Denote by $V^{\cal S}_{\lambda}$ the only subspace from $\mathcal V_\lambda(\mathcal S)$ that depends on the parameter $\lambda$ and denote by $a$ the element from $\cal S$ that corresponds to $V^{\cal S}_{\lambda}$. Deleting $V^{\cal S}_{\lambda}$ from $\mathcal V_\lambda(\mathcal S)$,  we obtain the subspace representation $\mathcal V(\mathcal S_a)=(V^{\cal S};V^{\cal S}_i)_{i \in \mathcal S_a}$ of primitive poset $\mathcal S_a:={\cal S} \setminus \{ a \}$.

\begin{proposition} \label{NonEquiv}
$\mathcal V_\lambda(\mathcal S)$ is not equivalent to $\mathcal V_\mu(\mathcal S)$ if $\lambda \neq \mu$ for each critical poset $\cal S$. All subspace representation $\mathcal V(\mathcal S_a)$ are schurian.
\end{proposition}

\begin{proof} This proposition is proved by straightforward computations.
\end{proof}

Let $\cal S$ be a critical poset. The poset $\mathcal S_a$ is primitive and does not contain any of the critical posets, its subspace representation $\mathcal V(\mathcal S_a)$ is schurian.  By \cite[Proposition 3.1]{GrushYus}, there exists a weight which we denote by $\chi^a$, such that  $\mathcal V(\mathcal S_a)$ is $\chi^a$-stable. Write
\[
  R:=\min \Big\{\dim M-\sum_{i\in \mathcal S_a} \chi^a_i \dim (V_i^{\mathcal S} \cap M) \, \Big | \,  M \mbox{ is a proper subspace of } V^{\mathcal S} \Big\}.
\]
The subspace representation $\mathcal V({\mathcal S_a})$ is $\chi^a$-stable, hence $R>0$. Let $\varepsilon$ be such that
$R>\varepsilon>0$. Write $T:=1+(R-\varepsilon)\dim V_\lambda^{\mathcal S}(\dim V^\mathcal S)^{-1}$ and 
\begin{equation*}
\chi^{\mathcal S}=(\chi_i^{\mathcal S})_{i\in \mathcal S}, \quad 
        \chi^{\mathcal S}_i:=\begin{cases} \chi^a_i\cdot T^{-1}, & \mbox{if} \ \ i \in \mathcal S_a,\\
         (R-\varepsilon)\cdot T^{-1}, &  \mbox{if} \ \ i=a. \end{cases}
\end{equation*}
\begin{proposition} \label{StabCrit}
The subspace representations $\mathcal V_{\lambda}(\mathcal S)$ are $\chi^{\mathcal S}$-stable for all $\lambda$ and $\mathcal S$.
\end{proposition}
\begin{proof}
Note that
\begin{align*}
	\sum_{i \in \mathcal S}\chi_i^{\mathcal S}\dim V_i ^{\mathcal S}&=
	T^{-1}\sum_{i \in \mathcal S_a}\chi_i^{a}\dim V_i ^{\mathcal S}+\chi^{\mathcal S}_a \dim V_\lambda^{\mathcal S}\\
	&=T^{-1}\dim V^{\mathcal S}+(1-T^{-1})\dim V^{\mathcal S}=\dim V^{\mathcal S}.
\end{align*}
Let $M$ be any proper subspace of $V^{\mathcal S}$. Then
   \begin{equation*}
    \begin{split}
       \sum_{i \in \mathcal S}\chi_i^{\mathcal S}\dim
        (V_i^{\mathcal S} \cap M)&=T^{-1}\sum_{i \in \mathcal S_a}\chi_i^a \dim
        (V_i ^{\mathcal S} \cap M)+\chi_a^{\mathcal S} \dim(V_\lambda^{\mathcal S}\cap M) \\
        &\leq T^{-1}\left (\dim M(1-R)+ (R-\varepsilon) \dim(V_\lambda^{\mathcal S}\cap M) \right)\\
        &<T^{-1}\dim M (1-\varepsilon)<\dim M.
    \end{split}
    \end{equation*}
    Hence $\mathcal V_\lambda(\mathcal S)$ is $\chi^{\mathcal S}$-stable.
\end{proof}

\begin{proposition}\label{crInf}
Critical posets are unitary representation-infinite.
\end{proposition}
\begin{proof}
By Proposition \ref{NonEquiv} and Proposition \ref{StabCrit}, any critical poset $\mathcal S$ has an infinite number of nonequivalent $\chi^{\mathcal S}$-stable subspace representations. By Theorem \ref{stability-un}, $\mathcal S$ has an infinite number of indecomposable unitarily nonequivalent  $\chi^{\mathcal S}$-representations.
\end{proof}

\begin{proposition} \label{finCri}
If a poset $\mathcal S$ contains a critical poset (as a full
subposet), then $\mathcal S$ is unitary representation-infinite.
\end{proposition}

\begin{proof}
Suppose that $\mathcal S$ contains a critical poset $\mathcal S_c$. By Proposition \ref{crInf}, there exists a weight $\mathcal \chi^{c}$ such that $\mathcal S_c$ has an infinite number of indecomposable unitarily nonequivalent $\chi^c$-representations of dimension $d$. Define the following subset of $\mathcal S$:
\begin{equation*}
\begin{split}
  \mathcal S_{\max}&:=\{a \in \mathcal S \ |\ b \prec a \ \mbox{for some} \ b \in \mathcal S_c\ \}.\\
\end{split}
\end{equation*}
   For each $\mathcal \chi^{c}$-representation $\mathcal U=(U;U_i)_{i\in\mathcal S_c}$,
   define the unitary representation $\mathcal U^\prime=(U;U_i^\prime)_{i\in\mathcal S}$ of $\mathcal S$ as follows:
\begin{equation*}
        U_i^\prime:=\begin{cases}
                         0, &  \mbox{if} \ \ i \notin \mathcal S_{\max}\cup \mathcal S_{c}, \\
                         U_i, &  \mbox{if} \ \ i \in \mathcal S_c, \\
                         U, &  \mbox{if} \ \ i \in \mathcal S_{\max}.
                       \end{cases}
\end{equation*}
It is easy to check that $\mathcal U^\prime$ is $\chi^\prime$-representation, in which $\chi^\prime=(\chi_i^\prime)_{i \in \mathcal S}$ is defined  by
   \begin{equation*}
        \chi_i^\prime:=\begin{cases}
                               \chi^c_i\cdot(1+|\mathcal S_{\text {max}}|)^{-1}, & \mbox{if} \ \   i \in \mathcal S^c, \\
                                (1+|\mathcal S_{\text{max}}|)^{-1}, &  \mbox{otherwise}.
                       \end{cases}
   \end{equation*}
   Hence $\mathcal S$ is unitary representation-infinite.
\end{proof}

The implication ${\rm(ii)}\Rightarrow {\rm(iii)}$ follows from Proposition \ref{finCri}. This finishes the proof of Theorem \ref{maintheoremFin}.

\begin{remark} \label{RemFinale}
\noindent Define the following weights:
 \begin{align*}
    \chi^{(1,1,1,1)}&:=\frac 12(1,1,1,1),\\
    \chi^{(2,2,2)}&:=\frac 13(1,1,1,1,1,1),\\
    \chi^{(1,3,3)}&:=\frac 14(2,1,1,1,1,1,1),\\
    \chi^{(1,2,5)}&:=\frac 16(3,2,2,1,1,1,1,1),\\
    \chi^{(N,4)}&:=\frac 15(2,1,1,2,1,1,1,1).
 \end{align*}
Each weight $\chi^{\mathcal S}$ obtained from the minimal imaginary root of the quadratic form related to a critical poset $\mathcal S$.  We checked (describing all possible subdimension vectors) that the representations $\mathcal V_\lambda(\mathcal S)$ are $\chi^{\mathcal S}$-stable for any $\lambda \in \mathbb C \setminus \{0,1\}$. Hence they give rise to an infinite family of nonequivalent
$\chi^{\mathcal S}$-representations.
For primitive $\mathcal S$  one can obtain the precise description of projections for such representations using the results from \cite{AlbeverioOstrSamoilenko, KruglyakLivinskyi,OstrovskyiSamoilenko}. The description in the case $\mathcal S=(N,4)$ is unknown.
\end{remark}

\medskip

\textbf{Acknowledgments.} We would like to thank S.A. Kruglyak,
V.L. Ostrovskii, and L.B. Turowska for
useful discussions and advises and also V.V. Sergeichuk for numerous helpful remarks.
The second author was partially supported by DFG grant SCHM1009/4-1 and Fapesp grant 2010/15781-0.

%\bibliographystyle{plain}
%\bibliography{bibliography}

\begin{thebibliography}{10}
%\end{thebibliography}


\bibitem{AlbeverioOstrSamoilenko}
S.~Albeverio, V.~Ostrovskyi, Yu. Samoilenko,
\newblock On functions of graphs and representations of a certain class of
  $*$-algebras,
\newblock {J. Algebra} 308 (2) (2007) 567--582.


\bibitem{GabRoi}
P.~Gabriel, A.~Roiter,
\newblock Representations of finite-dimensional algebras,
Encyclopaedia of mathematical sciences 73, Springer, 1997.


\bibitem{GrushYus}
R.~Grushevoy, K.~Yusenko,
\newblock On the unitarization of linear representations of primitive partially ordered sets,
\newblock {Operator Theory: Advances and Applications} 190 (2) (2009) 279--294.


\bibitem{hal}
P.R.~Halmos, Two subspaces, {Trans. Amer. Math. Soc.} 144 (1969) 381--389.


\bibitem{Hu2004}
Yi~Hu,
\newblock Stable configurations of linear subspaces and quotient coherent
  sheaves,
\newblock {Q. J. Pure Appl. Math.} 1 (1) (2005) 127--164.


\bibitem{Kleiner1}
M.M.~Kleiner,
\newblock Partially ordered sets of finite type,
\newblock {J. Soviet Math.} 3 (1975)  607--615.

\bibitem{Klyachko}
A.A.~Klyachko,
\newblock Stable bundles, representation theory and hermitian operators,
\newblock {Selecta Math.} 4 (1988) 419--445.

\bibitem{KruglyakLivinskyi}
S.A.~Kruglyak, I.V.~Livinskyi,
\newblock Regular orthoscalar representations of the extended Dynkin graph
  $\widetilde E_8$ and $*$-algebra associated with it,
\newblock {Ukrainian Math. J.} 62 (8) (2010) 1213--1233.

\bibitem{KruglyaNazarovaRoiter}
S.A.~Kruglyak, L.A.~Nazarova, A.V.~Roiter,
\newblock Orthoscalar quiver representations corresponding to extended Dynkin
  graphs in the category of Hilbert spaces,
\newblock {Funct. Anal. Appl.} 44 (2) (2010) 125--138.

\bibitem{KruglyaSamoilenko2}
S.A.~Kruglyak, Yu.S.~Samoilenko,
\newblock On the complexity of description of representations of $*$-algebras
  generated by idempotents,
\newblock {Proc. AMS}
  128 (6) (2000) 1655--1664.


\bibitem{NazarovaRoiter}
L.A.~Nazarova, A.V.~Roiter,
\newblock Representations of partially ordered sets,
\newblock Zap. Nauchn. Sem. Leningrad. Otdel. Mat. Inst.
Steklov. (LOMI) 28, 1972, 5--31 (in Russian);  English
version: J. Soviet Math. 3 (no. 5) (1975) 585--606.


\bibitem{OstrovskyiSamoilenko}
V.L.~Ostrovskyi, Yu.S.~Samoilenko,
\newblock On spectral theorems for families of linearly connected selfadjoint
  operators with prescribed spectra associated with extended Dynkin graphs,
\newblock {Ukrainian Math. J.} 58 (11) (2006) 1768--1785.

\bibitem{Sergeichuk}
V.V.~Sergeichuk,
\newblock Unitary and Euclidean representations of a quiver,
\newblock {Linear Algebra Appl.} 278 (1998) 37--62.

\bibitem{Simson1992}
D. Simson,
\newblock {Linear representations of partially ordered sets and vector
  space categories,} Algebra, Logic and Applications 4, Gordon \& Breach Science Publishers, 1992.


\bibitem{Shapiro}
H.~Shapiro,
\newblock A survey of canonical forms and invariants for unitary similarity,
\newblock {Linear Algebra Appl.} 147 (1991) 101--167.

\bibitem{Totaro1994}
B.~Totaro,
\newblock Tensor products of semistables are semistable,
\newblock {Geometry and Analysis on Complex Manifolds, World Series, River
  Edge}, 1994, 242--250.


\bibitem{Wu}

P.Y.~Wu,
\newblock Sums of idempotent matrix,
\newblock Linear Algebra Appl. 142 (1990) 43--54.




\end{thebibliography}

\end{document}